\documentclass[runningheads]{llncs}
\usepackage{amssymb}
\usepackage{amsmath}
\usepackage{graphicx}
\usepackage{graphics}
\usepackage{adjustbox}
\usepackage{lineno}
\usepackage{color}
\usepackage{algpseudocode}
\usepackage{algorithm}

\algnewcommand\algorithmicinput{\textbf{Input:}}
\algnewcommand\Input{\item[\algorithmicinput]}
\algnewcommand\algorithmicoutput{\textbf{Output:}}
\algnewcommand\Output{\item[\algorithmicoutput]}

\DeclareMathOperator*{\argmin}{arg\,min}

\begin{document}
			
\title{A study on a feedforward neural network to solve partial differential equations in \\ hyperbolic-transport problems\thanks{Supported by S\~ao Paulo Research Foundation (FAPESP), National Council for Scientific and Technological Development, Brazil (CNPq), and PETROBRAS - Brazil.}}
\titlerunning{A study on a feedforward neural network to solve PDEs}

\author{Eduardo Abreu\inst{1}\orcidID{0000-0003-1979-3082} \and Joao B. Florindo\inst{1}\orcidID{0000-0002-0071-0227}}

\authorrunning{Abreu and Florindo}

\institute{Institute of Mathematics, Statistics and Scientific Computing - University of Campinas,
	Rua S\'{e}rgio Buarque de Holanda, 651, Cidade Universit\'{a}ria "Zeferino Vaz" - Distr. Bar\~{a}o Geraldo, CEP 13083-859, Campinas, SP, Brasil
	%\email{lncs@springer.com}\\
	\url{http://www.ime.unicamp.br}\\
	\email{\{eabreu,florindo\}@unicamp.br}}

\maketitle

\begin{abstract}
In this work we present an application of modern deep lear\-ning methodologies to the numerical solution of partial differential equations in transport models. More specifically, we employ a supervised deep neural network that takes into account the equation and initial conditions of the model. We apply it to the Riemann problems over the inviscid nonlinear Burger's equation, whose solutions might develop discontinuity (shock wave) and rarefaction, as well as to the classical one-dimensional Buckley-Leverett two-phase problem. The Buckley-Leverett case is slightly more complex and interesting because it has a non-convex flux function with one inflection point. Our results suggest that a relatively simple deep learning model was capable of achieving promising results in such challenging tasks, providing numerical approximation of entropy solutions with very good precision and consistent to classical as well as to recently novel numerical methods in these particular scenarios.
\end{abstract}

\section{Introduction}

Deep learning techniques have been applied to a variety of problems in science during the last years, with numerous examples in image recognition \cite{HZRS16}, natural language processing \cite{YHPC18}, self driving cars \cite{CSKX15}, virtual assistants \cite{KB18}, healthcare \cite{LKBSC17}, and many others. More recently, we have seen a growing interest on applying those techniques to the most challenging problems in mathematics and the solution of differential equations, especially partial differential equations (PDE), is a canonical example of such task \cite{RPK19}.

Despite the success of recent learning-based approaches to solve PDEs in rela\-tively ``well-behaved'' configurations, we still have points in these methodologies and applications that deserve more profound discussion, both in theoretical and practical terms. One of such points is that many of these models are based on complex structures of neural networks, sometimes comprising a large number of layers, recurrences, and other ``ad-hoc'' mechanisms that make them difficult to be trained and interpreted. Furthermore, we have seen little discussion on more challenging problems, like those involving discontinuities and ``shock'' solutions numerical approximation of entropy solutions in hyperbolic-transport problems, see, e.g., \cite{ABO2020,GCJG19,CD2016,HKRS19,OAO63,DSLS19,EAJP19,ADGP20} and references cited therein.

This is the motivation for the study accomplished in this work, where we investigate a simple feed-forward architecture, based on the physics-informed model proposed in \cite{RPK19}, to complex problems involving PDEs in transport mo\-dels. More specifically, we analyze the numerical solutions of four initial-value problems: three problems on the inviscid nonlinear Burgers PDE (involving shock wave and smooth/rarefaction fan for distinct initial conditions) and on the one-dimensional Buckley-Leverett equation for two-phase configuration, which is a litle rather more complex and interesting because it has a non-convex flux function with one inflection point. The neural network consists of 9 stacked layers with $\tanh$ activation and geared towards minimizing the approximation error both for the initial values and for values of the PDE functional calculated by automatic differentiation.

The achieved results are promising. We managed to obtain an average quadra\-tic error of 0.0005, 0.0018, 0.0001, and 0.0021, respectively, for the rarefaction, shock, smooth, and Buckley-Leverett problems. Such results are pretty inte\-resting if we consider the low complexity of the neural model and the challenge involved in these discontinuous cases. It also strongly suggests more in-depth studies on deep learning models that account for the underlying equation. They seem to be a quite promising line to be explored for challenging problems arising in physics, engineering, and many other areas.

\section{Hyperbolic problems in transport models}

Hyperbolic partial differential equations in transport models describe a wide range of wave-propagation and transport phenomena arising 
from scientific and industrial engineering area. This is a fundamental research that is in active progress since it involves complex multiphysics and advanced simulations due to a lack of general mathematical theory for closed-analytical solutions (see, e.g., \cite{HKRS19,CLHK20,BCS2020,ABO2020,GCJG19,DSLS19,CD2016} and references cited therein). A basic fact of nonlinear hyperbolic transport problems is the possible loss of regularity in their solutions, namely, even solutions which are initially smooth (i.e., initial datum) may become discontinuous within finite time (blow up in finite time) and after this singular time, nonlinear interaction of shocks and rarefaction waves will play an important role in the dynamics. For the sake of simplicity and without loss of generality,  we consider the scalar 1D Cauchy problem 
\begin{equation}
\frac{\partial u}{\partial t}+
\frac{\partial H(u)}{\partial x} = 0, 
\quad x \in \mathbb{R}, \quad t > 0, \quad \quad \quad u(x,0) = u_0(x),
\label{ivp2b}
\end{equation}
here $H \in C^{2}(\Omega)$, $H:\Omega \to \mathbb{R}$,
$u_{0}(x) \in L^{\infty}(\mathbb{R})$ and $u=u(x,t): \mathbb{R} 
\times \mathbb{R^+} \longrightarrow \ \Omega \subset 
\mathbb{R}.$ In several applications, the flux function $H(u)$ 
is smooth and with a finite number of inflection points, namely
$H(u)=\frac{u^2}{u^2+a(1-u)^2}$, $0 < a < 1$, such as is the 
case of classical Buckley-Leverett equation (it has a non-convex 
flux function with one inflection point) in petroleum engineering 
(e.g., \cite{ADGP20}). Another interesting model is the inviscid 
Burgers' equation where $H(u)=u^2/2$, which is used to model for 
example gas dynamics and traffic flow, and investigate the appearance 
of shock waves, especially in fluid mechanics, for nonlinear wave 
propagation (see, e.g., \cite{EAJP19}).

By using an argument in terms of traveling waves to capture the 
viscous profile at shocks, one can conclude that solutions of 
(\ref{ivp2b}) satisfy Oleinik's entropy condition (see e.g., 
\cite{OAO63}), which are limits of solutions 
$u^{\epsilon}(x,t) \rightarrow u(x,t)$, where $u(x,t)$ is
given by (\ref{ivp2b}) and $u^{\epsilon}(x,t)$ is given by
the augmented parabolic equation \cite{BKQ71}
\begin{equation}
\frac{\partial u^{\epsilon}}{\partial t}
+ \frac{\partial H(u^{\epsilon})}{\partial x} 
= \epsilon \frac{\partial^2 u^{\epsilon}}{\partial x^2},
\quad x \in \mathbb{R}, \quad t > 0, \quad \quad \quad u^{\epsilon}(x,0) = u^{\epsilon}_0(x),
\label{ivp2bEps}
\end{equation}
with $\epsilon > 0$ and the same initial data as in (\ref{ivp2b}).

Thus, in many situations is of importance to consider and study
both model problems (\ref{ivp2b}) and (\ref{ivp2bEps}), related
to hyperbolic problems in transport models. In this regard, a 
typical flux function $H(u)$ associated to fundamental prototype 
models (\ref{ivp2b}) and (\ref{ivp2bEps}) depends on the application 
under consideration, for instance, such as modeling slow erosion 
phenomena in granular flow, fluid mechanics, flow in porous media 
(see, e.g., \cite{EAJP19,ADGP20,ABO2020,HKRS19} and references 
cited therein). Moreover, it is noteworthy that in practice 
calibration of function $H(u)$ can be difficult to achieve due 
to unknown parameters and, thus, data assimilation (i.e., 
regression method) can be an efficient method of calibrating 
these flux function $H(u)$ mo\-dels \cite{JBKN19,Rudye1602614}. 
We intend to design a unified approach which combines both 
Partial Differential Equation (PDE) modeling and fine tuning 
machine learning techniques aiming as a fisrt step to an 
effective tool for advanced simulations related to hyperbolic 
problems in transport models.

\subsection{A benchmark numerical scheme for solving model (\ref{ivp2b})}

From the above discussion, and for comparison purposes, we 
will provide correct qualitative entropy approximation solutions 
for model problem (\ref{ivp2b}) by using two numerical schemes, 
namely, the conservative method 
\begin{equation}
	U_j^{n+1} = U_j^n - \frac{k}{h} \,
	\Big[ \, F(U_j^n,U_{j+1}^n) - F(U_{j-1}^n,U_{j}^n) \, \Big],
	\label{LFx}
\end{equation}
with the associated classical Lax-Friedrichs numerical flux
found elsewhere,
\begin{equation}
	F(U_j^n,U_{j+1}^n) = \frac{1}{2} \,
	\left[
		\frac{h}{k} \, (U_j^n - U_{j+1}^n) +
		\Big( \, H(U_{j+1}^n) + H(U_j^n) \, \Big)
	\right],
	\label{LFxflux}
\end{equation}
as well as by the the Lagrangian-Eulerian numerical flux (\cite{EAJP19,AMPR20})
\begin{equation}
	F(U_j^n,U_{j+1}^n) = {\frac{1}{4}}\,
	\left[
		\frac{h}{k}\,(U_j^n - U_{j+1}^n) + 
		2 \, \Big( \, H(U_{j+1}^n) + H(U_j^n) \, \Big)
	\right].
	\label{ELDicFlux}
\end{equation}
Here, both schemes (\ref{LFxflux}) and (\ref{ELDicFlux}) should
follow the stability Courant--Friedrichs--Lewy (CFL) condition
%%%
\begin{equation}
	\label{CFLcondLB}
	\max\limits_{j}\{ \, |H'(U_j^n)| \, \} \, \frac{k}{h} \,<\, \frac{1}{2},
\end{equation}
for all time steps $n$, where $k = \Delta t^n$ and 
$h=\Delta x$, $H'(U_j^n)$ is the partial derivative of 
$H$, namely $\displaystyle\frac{\partial H(u)}{\partial u}$
for all $U_j^n$ in the mesh grid.

\section{Related works}

We are interested in the study of a unified approach which 
combines both data-driven models (regression method by machine 
learning) and physics-based models (PDE modeling). Our approach 
is substantially distinct from the current trend of merely 
data-driven discovery type methods for recovery governing 
equations by using machine learning and artificial intelligence 
algorithms in a straightforward manner. We glimpse the use of 
novel methods, fine tuning machine learning algorithms and very 
fine mathematical and numerical analysis to improve comprehension 
of regression methods aiming to identify the potential and 
reliable prediction for advanced simulation for hyperbolic 
problems in transport models as well as the estimation of 
financial returns and economic benefits. In this regard, we 
mention the very interesting works 
\cite{RDQ20,GKO20,JBKN19,Rudye1602614,RPK19}, which 
introduced physics-informed neural networks that are trained 
to solve supervised learning tasks while respecting conceptual 
foundations of physics as well as proving data assimilation 
for summarizing data-driven discovery of partial differential 
equations. Related to the hyperbolic problems in transport 
models (\ref{ivp2b}) and (\ref{ivp2bEps}), we mention the 
very recent review paper \cite{BNK20}, which discusses machine 
learning for fluid mechanics, but highlighting that such 
approach could augment existing efforts for the study, 
modeling and control of fluid mechanics, keeping in mind 
the importance of honoring engineering principles and the 
governing equations, mathematical and physical foundations 
driven by unprecedented volumes of data from experiments 
and advanced simulations at multiple spatiotemporal scales. 
We also mention the work \cite{RBD20}, where the issue of 
{\it domain knowledge} is addressed as a prerequisite 
essential to gain explainability to enhance scientific 
consistency from machine learning and foundations of 
physics-based given in terms of mathematical equations 
and physical laws. However, we have seen much less 
discussion on more challenging PDE modeling problems, 
like those involving discontinuities and shock' solutions 
numerical approximation of entropy solutions in 
hyperbolic-transport problems, in which the issue of 
conservative numerical approximation of entropy solutions 
is crucial and mandatory 
\cite{ABO2020,GCJG19,CD2016,HKRS19,OAO63,DSLS19,EAJP19,ADGP20,AMPR20}.

\section{Proposed methodology}

The neural network employed here is based on that described in \cite{RPK19}. It follows a classical feed-forward architecture, with 9 hidden layers, each one with a hyperbolic tangent used as activation function.

The general problem solved here has the form
\begin{equation} u_t + \mathcal{N}(u) = 0, \qquad x \in \Omega, t \in [0,T],\end{equation}
where $\mathcal{N}(\cdot)$ is a non-linear operator and $u(x,t)$ is the desired solution. Unlike the methodology described in \cite{RPK19}, here we do not have an explicit boundary condition and the neural network is optimized only over the initial conditions of each problem.

We focus on four problems: the inviscid nonlinear Burgers equation 
\begin{equation} 
	u_t + \left( \frac{u^2}{2} \right)_x = 0, \qquad x\in[-10,10], \qquad t \in [0,8],
\end{equation}
with shock initial condition
\begin{equation} u(x,0) = 1, x < 0 \mbox{ and } u(x,0) = 0, x > 0, \end{equation}
discontinuous initial data (hereafter {\it rarefaction fan initial condition})
\begin{equation} u(x,0) = -1, x < 0 \mbox{ and } u(x,0) = 1, x > 0, \end{equation}
smooth initial condition
\begin{equation} u(x,0) = 0.5 + \sin(x), \end{equation}
and the two-phase Buckley-Leverett
\begin{equation} 
	\begin{split}
		&u_t + \left( \frac{u^2}{u^2+a(1-u)^2} \right)_x = 0, \qquad x\in[-8,8], \qquad t \in [0,8],\\
		&u(x,0) = 1, x < 0 \mbox{ and } u(x,0) = 0, x > 0.\\
	\end{split}
\end{equation}
In this problem we take $a=1$.

For the optimization of the neural network we should define $f$ as the left hand side of each PDE, i.e.,
\begin{equation} f := u_t + \mathcal{N}(u), \end{equation}
such that 
\begin{equation} \mathcal{N}(u) = \left( \frac{u^2}{2} \right)_x \end{equation}
in the inviscid Burgers and
\begin{equation} \mathcal{N}(u) = \left( \frac{u^2}{u^2+a(1-u)^2} \right)_x \end{equation}
in the Buckley-Leverett. Here we also have an important novelty which is the introduction of a derivative (w.r.t. $x$) in $\mathcal{N}(u)$, which was not present in \cite{RPK19}.

The function $f$ is responsible for capturing the physical structure (i.e, select the qualitatively correct entropy solution) of the problem and inputting that structure as a primary element of the machine learning problem. The neural network computes the expected solution $u(x,t)$ and its output and the derivatives present in the calculus of $f$ are obtained by automatic differentiation. 

Two quadratic loss functions are defined over $f$, $u$ and the initial condition:
\begin{equation} 
	\begin{split}
		& \mathcal{L}_f(u) = \frac{1}{N_f}\sum_{i=1}^{N_f}|f(x_f^i,t_f^i,)|^2,\\	
		& \mathcal{L}_u(u) = \frac{1}{N_u}\sum_{i=1}^{N_u}|u(x_u^i,t_u^i)-u^i|^2,\\
	\end{split}
\end{equation}
where $\{x_f^i,t_f^i\}_{i=1}^{N_f}$ correspond to collocation points over $f$, whereas $\{x_u^i,t_u^i,u^i\}_{i=1}^{N_u}$ correspond to the initial values at pre-defined points. 

Finally, the solution $u(x,t)$ is approximated by minimizing the sum of both objective functions at the same time, i.e.,
\begin{equation} u(x,t) \approx \argmin_{u} [\mathcal{L}_f(u) + \mathcal{L}_u(u)]. \end{equation}

\section{Results and Discussion}

In the following we present results for the solutions of the investigated problems obtained by the neural network model. We compare these solutions with two numerical schemes: Lagrangian-Eulerian and Lax-Friedrichs. These are very robust numerical methods with a solid mathematical basis. Here we use one scheme to validate the other. In fact, the solutions obtained by each scheme are very similar. For that reason, we opted for graphically showing curves only for the Lagrangian-Eulerian solution. However, we exhibit the errors of the proposed methodology both in comparison with Lagrangian-Eulerian (EEL) and Lax-Friedrichs (ELF). Here such error corresponds to the average quadratic error, i.e.,
\begin{equation}
	\begin{split}	
	& ELF(t) = \frac{\sum_{i=1}^{N_u}(u_{NN}(x^i,t)-u_{LF}(x^i,t))^2}{N_u}, \\
	& EEL(t) = \frac{\sum_{i=1}^{N_u}(u_{NN}(x^i,t)-u_{LE}(x^i,t))^2}{N_u},
	\end{split}
\end{equation}
where $u_{NN}$, $u_{LF}$, and $u_{LE}$ correspond, respectively, to the neural network, Lax-Friedrichs, and Lagrangian-Eulerian solutions. In our tests, we used $N_f=10^4$ unless otherwise stated, and $N_u=100$. For the numerical reference schemes we adopted CFL condition 0.4 for Lax-Friedrichs and 0.2 for Lagrangian-Eulerian. We also used $\Delta x = 0.01$.

\begin{figure}[!htpb]
	\centering
	\includegraphics[width=\textwidth]{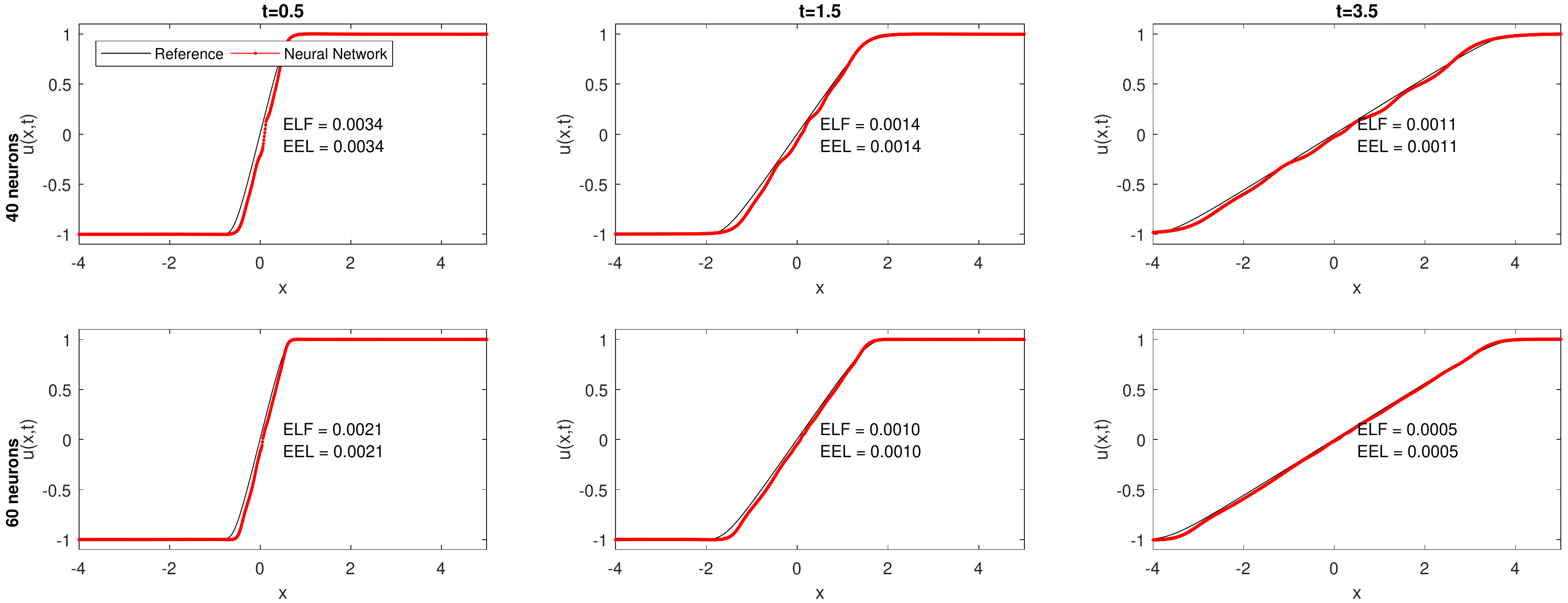}
	\caption{Burgers: Rarefaction.}
	\label{fig:rarefaction}
\end{figure}

For the rarefaction case, we observed that using $N_f=10^4$ collocation points was sufficient to provide good results. In this scenario, we also verified the number of neurons, testing 40 and 60 neurons. Figure \ref{fig:rarefaction} shows the obtained solution compared with reference and the respective errors. Interestingly, the error decreases when time increases, which is a consequence of the solution behavior, which becomes smoother (smaller slope) for larger times, showing good accuracy and evidence that we are computing the correct solution in our numerical simulation.

\begin{figure}[!htpb]
	\centering
	\includegraphics[width=\textwidth]{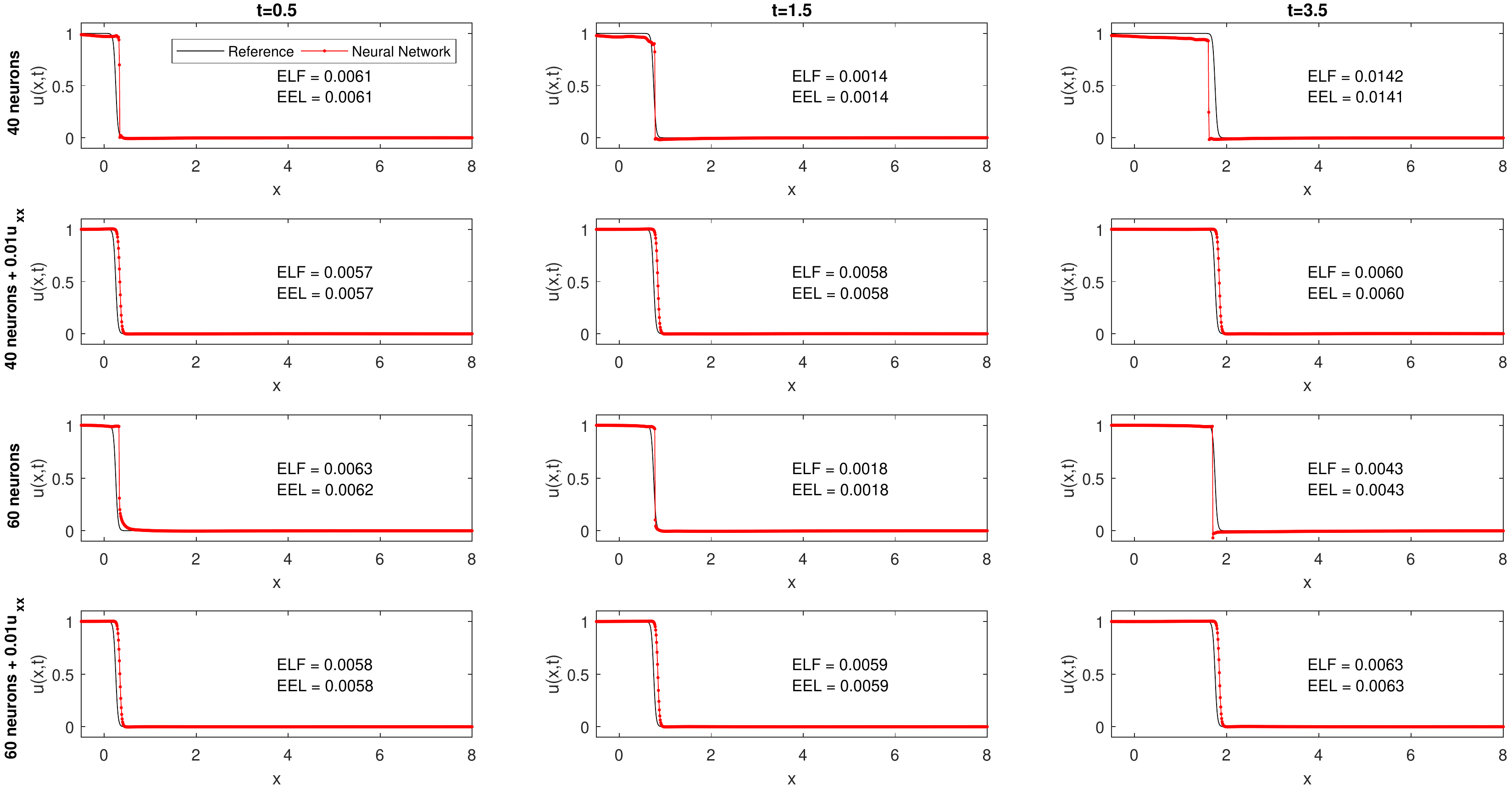}
	\caption{Burgers: Shock.}
	\label{fig:shock}
\end{figure}

Figure \ref{fig:shock} illustrates the performance of the neural network model for the inviscid Burgers equation with shock initial condition. Here we had to add a small viscous term ($0.01u_{xx}$) for better stabilization, but in view on the modeling problems (\ref{ivp2b}) and (\ref{ivp2bEps}). 
Here, such underlying viscous mechanism did not bring significant reduction in error, but the general structure of the obtained solution is better, attenuating spurious fluctuations around the discontinuities. It is crucial to mention at this point that numerical approximation of entropy solutions (with respect to the neural network) to hyperbolic-transport problems also require the notion of entropy-satisfying weak solution.  
It was also interesting to see that the addition of more neurons did not reduce the error for this initial condition. This is a typical example of overfitting caused by over-parameterization. An explanation for that is the relative simplicity of the initial condition, assuming only two possible values.

\begin{figure}[!htpb]
	\centering
	\includegraphics[width=\textwidth]{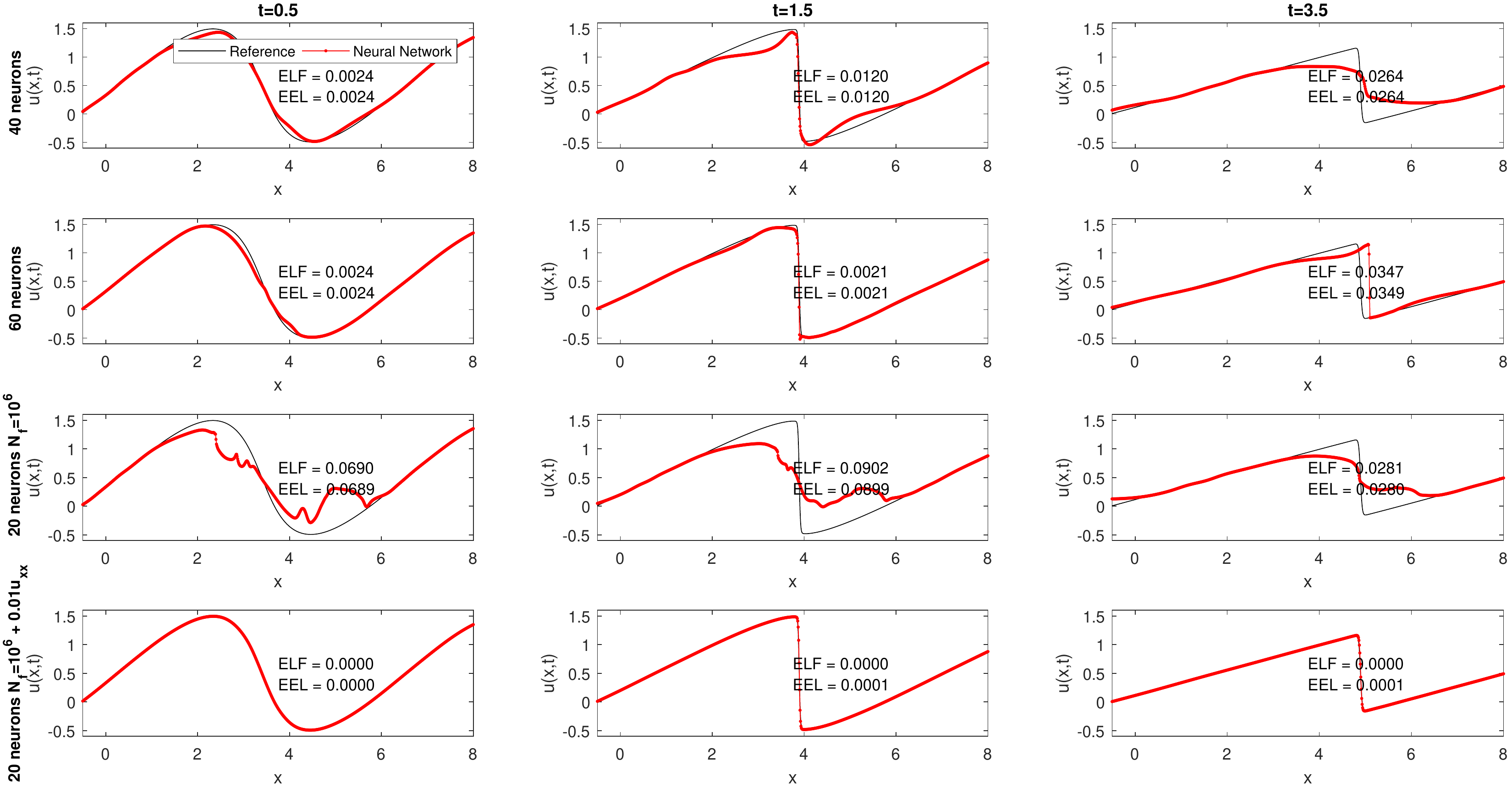}
	\caption{Burgers: Smooth.}
	\label{fig:smooth}
\end{figure}

Figure \ref{fig:smooth} depicts the solutions for the smooth initial condition in the inviscid Burgers equation. Here, unlike the previous case, increasing the number of neurons actually reduced the error. And this was expected considering that now both initial condition and solution are more complex. Nevertheless, we identified that tuning only number of neurons was not enough to achieve satisfactory solution in this situation. Therefore we also tuned the parameter $N_f$. In particular, we discovered that combining the same small viscous term used for the shock case with $N_f=10^6$ provided excellent results, with quite promising precision in comparison with our reference solutions.

\begin{figure}[!htpb]
	\centering
	\includegraphics[width=\textwidth]{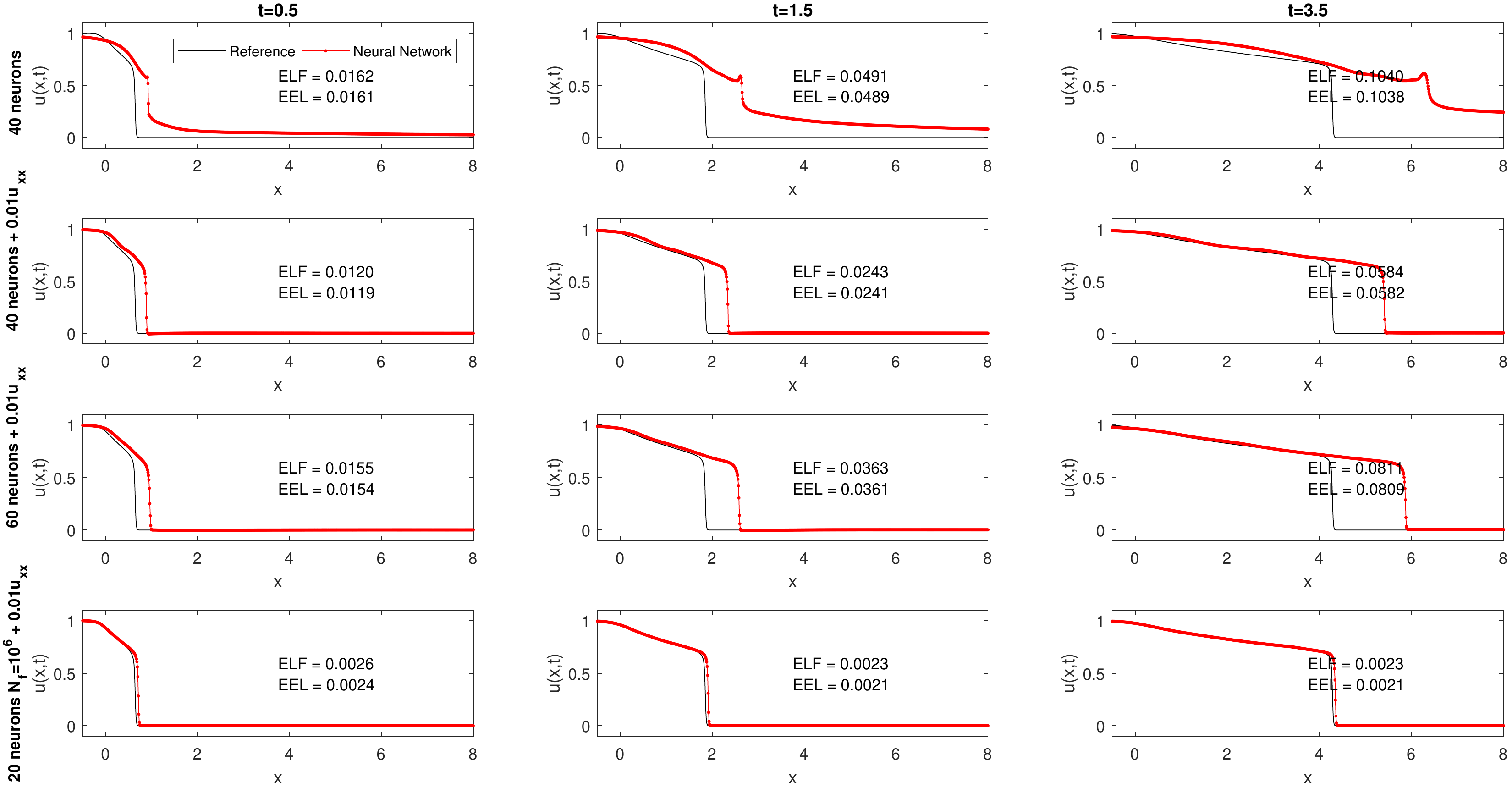}
	\caption{Buckley-Leverett: Rarefaction + Shock.}
	\label{fig:buckley}
\end{figure}

Another case characterized by solutions with more complex behavior is Buckley-Leverett with shock initial condition (Figure \ref{fig:buckley}). And, similar to what happened in the smooth case, again, the combination of $N_f=10^6$ with the small viscous term was more effective than any increase in the number of neurons. While the introduction of the small viscous term attenuated fluctuations in the solution when using 40 neurons, at the same time when using $N_f=10^4$, we observe that increasing the number of neurons causes an increase in the delay between the solution provided by the network and the reference. 

Generally speaking, the neural networks here studied were capable of achie\-ving promising results in challenging situations, involving different types of discontinuities and nonlinearities. Moreover, our results also bring important theoretical implications. In particular, the neural networks obtained results pretty close to those provided by entropic numerical schemes like Lagrangian-Eulerian and Lax-Friedrichs. Going beyond the analysis in terms of raw precision, these results give us evidences that our neural network model possess some type of entropic property, which from the viewpoint of a numerical method is a fundamental and desirable characteristic. 

\section{Conclusions}

This work presented an application of a feed-forward neural network to solve challenging hyperbolic problems in transport models. More specifically, we solve the inviscid Burgers equation with shock, smooth and rarefaction initial conditions, as well as the Buckley-Leverett equation with classical Riemann datum, which lead to the well-known solution that comprises (from left to right) a rare\-faction and a (discontinuous) shock wave. Our network was tuned according to each problem and interesting findings were observed. At first, our neural network model was capable of providing solutions pretty similar to those obtained by two numerical schemes here used as references: Lagrangian-Eulerian and Lax-Friedrichs. Besides, the general structure of the obtained solutions also behaved as expected, which considering the intrinsic challenge of these problems is a remarkable achievement. In fact, the investigated neural networks showed evidences of an entropic property, which is an important attribute of a numerical scheme, especially in problems like those investigated here. 

In summary, the obtained results share both practical and theoretical implications. In practical terms, the results confirm the potential of a relatively simple deep learning model in the solution of an intricate numerical problem. In theoretical terms, this also opens an avenue for formal as well as rigorous studies on these networks as mathematically valid and effective numerical methods. 

\section*{Acknowledgements}

J. B. Florindo gratefully acknowledges the financial support of S\~ao Paulo Research Foundation (FAPESP) (Grant \#2016/16060-0) and from National Council for Scientific and Technological Development, Brazil (CNPq) (Grants \#301480\\/2016-8 and \#423292/2018-8). E. Abreu gratefully acknowledges the financial support of S\~ao Paulo Research Foundation (FAPESP) (Grant \#2019/20991-8), from National Council for Scientific and Technological Development - Brazil (CNPq) (Grant \#2 306385/2019-8) and PETROBRAS - Brazil (Grant \#2015/\\00398-0). E. Abreu and J. B. Florindo also gratefully acknowledge the financial support of Red Iberoamericana de Investigadores en Matemáticas Aplicadas a Datos (MathData).

%\bibliographystyle{splncs04}
%\bibliography{MLPDE}

\end{document}